\documentclass[11pt]{article}
\usepackage{amssymb,amsfonts}
\begin{document}

\title{Note on the Painlev\'e V tau-functions}
\author{Yu.\,P.\,Bibilo, R.\,R.\,Gontsov}

\date{}

\maketitle

\begin{abstract}
We study some properties of tau-functions of an isomonodromic deformation leading to the fifth Painlev\'e
equation. In particular, here is given an elementary proof of Miwa's formula for the logarithmic differential
of a tau-function.
\end{abstract}

\section{Introduction}

This work is an addition to the article \cite{BG}, where we studied some
properties of the Malgrange isomonodromic deformation of a linear
differential $(2\times2)$-system defined on the Riemann sphere
$\overline{\mathbb C}$ and having at most two irregular
singularities of Poincar\'e rank one. Here we consider a
particular case of such a system:
\begin{eqnarray}\label{22syst}
\frac{dy}{dz}=\biggl(\frac{B^0_0}z+\frac{B^0_1}{z-t_0}
+\Bigl(\begin{array}{cc} 1 & 0 \\ 0 & 0
       \end{array}\Bigr)\biggr)y,
\end{eqnarray}
where $y(z)\in{\mathbb C}^2$ and $B^0_0$, $B^0_1$ are $(2\times
2)$-matrices. This system has two Fuchsian singular points $z=0$,
$z=t_0\in{\mathbb C}\setminus\{0\}$ and one {\it non-resonant} irregular
singularity $z=\infty$ of Poincar\'e rank one.

As known (see \cite[\S\S10,11]{Wa} or \cite[\S21]{IY}), in a
neighbourhood of (non-resonant) irregular singularity $z=\infty$
the system (\ref{22syst}) is formally equivalent to a system 
$$
\frac{d\tilde y}{dz}=\left(\Bigl(\begin{array}{cc} 1 & 0 \\ 0 & 0
                                 \end{array}\Bigr)+\frac Az\right)\tilde y,
$$
where $A$ is a diagonal $(2\times2)$-matrix. This means that there is an invertible 
matrix formal Taylor series $\widehat F(z)$ in $1/z$ begining with the identity 
matrix such that these two systems are connected by means of the transformation 
$y=\widehat F(z)\tilde y$ (and such a series $\widehat F$ is unique). Thus, the system (\ref{22syst})
possesses a uniquely determined formal fundamental matrix $\widehat Y(z)$ of the form
$$
\widehat Y(z)=\widehat F(z)\,z^A\,e^{{\rm diag}\,(z,0)}.
$$
According to Sibuya's sectorial normalization theorem (see \cite[Th. 21.13, Prop. 21.17]{IY}), a punctured 
neighbourhood of the point $z=\infty$ is covered by two sectors
$$
S_1=\Bigl\{\frac{\pi}2-\varepsilon<\arg z<\frac{3\pi}2+\varepsilon,\quad|z|>R\Bigr\}, \quad
S_2=\left\{-\frac{\pi}2-\varepsilon<\arg z<\frac{\pi}2+\varepsilon,\quad|z|>R\right\},
$$
with a sufficiently small $\varepsilon>0$ and sufficiently large $R>0$, such that in each $S_i$ there exists
a {\it unique} actual fundamental matrix 
$$
Y_i(z)=F_i(z)\,z^A\,e^{{\rm diag}\,(z,0)}
$$
of the system (\ref{22syst}) whose factor $F_i(z)$ has the asymptotic expansion $\widehat F(z)$ in $S_i$.
 
The intersection $S_1\cap S_2$ is a union of two sectors $\Sigma_1$, $\Sigma_2$,
$$
\Sigma_1=\Bigl\{\frac{\pi}2-\varepsilon<\arg z<\frac{\pi}2+\varepsilon,\quad|z|>R\Bigr\}, \quad
\Sigma_2=\Bigl\{\frac{3\pi}2-\varepsilon<\arg z<\frac{3\pi}2+\varepsilon,\quad|z|>R\Bigr\}.
$$
In the sector $\Sigma_1$ the fundamental matrices $Y_1$ and $Y_2$ necessarily differ by a constant invertible
matrix:
$$
Y_2(z)=Y_1(z)C_1, \qquad C_1\in{\rm GL}(2,{\mathbb C}), \quad z\in\Sigma_1.
$$
In the sector $\Sigma_2$, one similarly has 
$$
Y_1(z)=Y_2(z)C_2, \qquad C_2\in{\rm GL}(2,{\mathbb C}), \quad z\in\Sigma_2.
$$
The matrices $C_1, C_2$ are called (Sibuya's) Stokes matrices of the system (\ref{22syst}) at the non-resonant
irregular singular point $z=\infty$ (more precisely, the second Stokes matrix is $C_2\,e^{-2\pi iA}$).

Further we will focus on isomonodromic deformations of the system
(\ref{22syst}) which are closely related to the fifth Painlev\'e equation. 
According to M.\,Jimbo \cite{Ji} such an isomonodromic deformation is a family 
\begin{eqnarray}\label{isomfam}
\frac{dy}{dz}=\left(\frac{B_0(t)}z+\frac{B_1(t)}{z-t}+
\Bigl(\begin{array}{cc} 1 & 0 \\ 0 & 0
      \end{array}\Bigr)\right)y, \qquad B_{0,1}(t_0)=B^0_{0,1},
\end{eqnarray}
of differential systems holomorphically depending on the parameter $t\in D(t_0)$, 
where $D(t_0)\subset\mathbb C$ is a neighbourhood of the point $t_0$ and the matrix functions $B_0(t)$,
$B_1(t)$ are determined by the integrability condition $d\Omega=\Omega\wedge\Omega$
for the matrix meromorphic differential 1-form
$$
\Omega=\left(\frac{B_0(t)}z+\Bigl(\begin{array}{cc} 1 & 0 \\ 0 & 0
                                  \end{array}\Bigr)\right)dz+
			 \frac{B_1(t)}{z-t}\,d(z-t) 
$$
on the space $\overline{\mathbb C}\times D(t_0)$. The isomonodromy of the family (\ref{isomfam}) means that
a part of its {\it monodromy data} is independent of $t$. This part consists of the monodromy matrices corresponding
to some fundamental matrix $Y(z,t)$ of (\ref{isomfam}), Stokes matrices at the infinity and the connection matrix
between $Y$ and $Y_1$.

Assuming the eigenvalues of $B_0(t)=(b_0^{ij}(t))$, $B_1(t)=(b_1^{ij}(t))$ to be 
$\pm\frac12\,\theta_0$, $\pm\frac12\,\theta_1\not\in\frac12\,\mathbb Z$ (they do not depend on $t$) and 
$B_0(t)+B_1(t)$ equivalent to ${\rm diag}(-\frac12\,\theta_{\infty},\frac12\,\theta_{\infty})$, one 
has the function
$$
u(t)=\frac{b_1^{12}(t)\Bigl(b_0^{11}(t)+\frac12\,\theta_0\Bigr)}
{b_0^{12}(t)\Bigl(b_1^{11}(t)+\frac12\,\theta_1\Bigr)}
$$
to satisfy the fifth Painlev\'e equation
$$
\frac{d^2u}{dt^2}=\left(\frac1{2u}+\frac1{u-1}\right)\left(\frac{du}{dt}\right)^2-\frac1t\frac{du}{dt}+
\frac{(u-1)^2}{t^2}\left(\alpha u+\frac{\beta}u\right)+
\frac{\gamma u}t+\frac{\delta u(u+1)}{u-1},
$$
where
$$
\alpha=\frac18(\theta_0-\theta_1+\theta_{\infty})^2,\quad
\beta=-\frac18(\theta_0-\theta_1-\theta_{\infty})^2,\quad
\gamma=1-\theta_0-\theta_1, \quad \delta=-\frac12.
$$

There is also a geometric approach of B.\,Malgrange \cite{Ma} to the 
isomonodromic deformation of the system (\ref{22syst}) we referred 
to in \cite{BG}, but here we do not immerse in details of that approach
in view of the sufficiency of an analytic language for the present note.
 
According to the Miwa--Malgrange--Helminck--Palmer theorem (see \cite{Mi},  
\cite[\S3]{Ma}, \cite{He} or \cite[\S3]{Pa}) the matrix functions $B_0(t)$, $B_1(t)$ 
holomorphic in $D(t_0)$ can be extended meromorphically to the 
universal cover $\mathcal{D}\cong\mathbb C$ of the set ${\mathbb C}\setminus\{0\}$ 
of locations of the pole $t_0$. The set $\Theta\subset\cal D$ of the poles of 
the extended matrix functions $B_0$, $B_1$ (which may be empty) is usually called 
the {\it Malgrange $\Theta$-divisor} of the family (\ref{isomfam}). For $t^*\in\Theta$,
a {\it local $\tau$-function} of (\ref{isomfam}) is any holomorphic in $D(t^*)$ function 
$\tau^*$ such that $\tau^*(t^*)=0$. There exists a function $\tau$ (called a global 
$\tau$-{\it function} of the isomonodromic deformation (\ref{isomfam}) or of the fifth
Painlev\'e equation) holomorphic on the whole space $\cal D$ whose zero set coincides 
with $\Theta$. Thus, the set $\Theta$ has no limit points in $\cal D$. As follows from 
the results of \cite{BG}, all the zeros of this $\tau$-function are simple (at least in 
the case of the irreducible monodromy of (\ref{isomfam})). Here we give an elementary proof 
of Palmer's theorem \cite{Pa} concerning a $\tau$-function of the isomonodromic deformation 
(\ref{isomfam}).
\medskip

{\bf Theorem 1.} {\it A $\tau$-function of the fifth Painlev\'e equation satisfies the equality
$$
d\ln\tau(t)=\frac12{\rm res}_{z=t}{\rm tr}\Bigl(B(z,t)\Bigr)^2dt,
$$
where $B(z,t)$ denotes the coefficient matrix of the family $(\ref{isomfam})$.}
\medskip

The above formula is also referred to as a definition of a $\tau$-function (see \cite{JM}, \cite{Ji}). 
Having this definition T.\,Miwa \cite{Mi} has proved the analyticity of such a $\tau$-function on 
$\cal D$ (then its zeros are {\it a priori} included in $\Theta$).  

Another definition of a $\tau$-function of the Painlev\'e V (and of the other Painlev\'e equations)
comes from a Hamiltonian form of this equation. As K.\,Okamoto has shown \cite{Ok}, there is a function
$\tau$ holomorphic on $\cal D$ whose logarithmic derivative is equal to a Hamiltonian along a solution
(this $\tau$-function depends on a solution, its zeros are included in the set of poles of a solution).

\section{Proof of Theorem 1}

Consider a point $t^*\in\Theta$ of the Malgrange $\Theta$-divisor of (\ref{isomfam}). Assuming that $t^*\ne-1$ we make the transformation 
$\xi=1/(z+1)$ of the independent variable\footnote{If $t^*=-1$ one makes a transformation $\xi=1/(z-c)$, $c\ne-1$.} and come from 
(\ref{isomfam}) to the isomonodromic family 
\begin{eqnarray}\label{isomfam2}
\frac{d\tilde y}{d\xi}=\left(\frac{B_0(t(s))}{\xi-1}+\frac{B_1(t(s))}{\xi-s}-
\frac1{\xi^2}\Bigl(\begin{array}{cc} 1 & 0 \\ 0 & 0
                   \end{array}\Bigr)-\frac{B_0(t(s))+B_1(t(s))}{\xi}\right)\tilde y, \quad\tilde y(\xi)=y(z(\xi)), 
\end{eqnarray}
depending on the parameter 
$$
s=\frac1{t+1}\quad\Bigl(\Longrightarrow t=\frac{1-s}s\Bigr),
$$ 
with the Fuchsian singular points $\xi=1$, $\xi=s$ and non-resonant irregular singularity $\xi=0$ of Poincar\'e rank 1. 
The infinity is a non-singular point of (\ref{isomfam2}), therefore this family is a particular case of those considered in the paper \cite{BG}. 
Using some properties of such families obtained in that paper we will prove Miwa's formula for a local $\tau$-function of (\ref{isomfam2}) in
a neighbourhood of the point $s^*=1/(t^*+1)\in s(\Theta)$. Namely, the following statement holds whose proof is presented in the next section.
\medskip

{\bf Lemma 1.} {\it There is a local $\tau$-function $\tilde\tau(s)$ of the family $(\ref{isomfam2})$ near $s^*$ such that
$$
d\ln\tilde\tau(s)=\frac12{\rm res}_{\xi=s}{\rm tr}\Bigl(\widetilde B(\xi,s)\Bigr)^2ds,
$$
where $\widetilde B(\xi,s)$ is the coefficient matrix of this family.} 
\medskip

Having Lemma 1 we conclude for the local $\tau$-function $\tilde\tau(s(t))$ of (\ref{isomfam}) near $t^*$ that
\begin{eqnarray*}
d\ln\tilde\tau(s(t))=\frac12{\rm res}_{\xi=s}{\rm tr}\Bigl(\widetilde B(\xi,s)\Bigr)^2(-s^2)dt.
\end{eqnarray*}
To connect ${\rm res}_{\xi=s}{\rm tr}\Bigl(\widetilde B(\xi,s)\Bigr)^2(-s^2)$ and ${\rm res}_{z=t}{\rm tr}\Bigl(B(z,t)\Bigr)^2$, 
let us compute the both. Since 
$$
{\rm res}_{z=t}\Bigl(B(z,t)\Bigr)^2=\frac{B_0B_1+B_1B_0}t+B_1\Bigl(\begin{array}{cc} 1 & 0 \\ 0 & 0 \end{array}\Bigr)+
\Bigl(\begin{array}{cc} 1 & 0 \\ 0 & 0 \end{array}\Bigr)B_1,
$$
one has
$$
\frac12{\rm res}_{z=t}{\rm tr}\Bigl(B(z,t)\Bigr)^2={\rm tr}\frac{B_0B_1}t+{\rm tr}\left(B_1\Bigl(\begin{array}{cc} 1 & 0 \\ 0 & 0 
                                                                                                 \end{array}\Bigr)\right).
$$
In a similar way,
$$
{\rm res}_{\xi=s}\Bigl(\widetilde B(\xi,s)\Bigr)^2=\frac{B_0B_1+B_1B_0}{s-1}-\frac{B_0B_1+B_1B_0+2B_1^2}s-\frac1{s^2}\left(
B_1\Bigl(\begin{array}{cc} 1 & 0 \\ 0 & 0 \end{array}\Bigr)+\Bigl(\begin{array}{cc} 1 & 0 \\ 0 & 0 \end{array}\Bigr)B_1\right),
$$
and
\begin{eqnarray*}
\frac12{\rm res}_{\xi=s}{\rm tr}\Bigl(\widetilde B(\xi,s)\Bigr)^2(-s^2)&=&\frac{s^2}{1-s}{\rm tr}(B_0B_1)+s\,{\rm tr}(B_0B_1)+
s\,{\rm tr}\,B_1^2+{\rm tr}\left(B_1\Bigl(\begin{array}{cc} 1 & 0 \\ 0 & 0 \end{array}\Bigr)\right)=\\
&=&\frac s{1-s}{\rm tr}(B_0B_1)+s\,{\rm tr}\,B_1^2+{\rm tr}\left(B_1\Bigl(\begin{array}{cc} 1 & 0 \\ 0 & 0 \end{array}\Bigr)\right)=\\
&=&{\rm tr}\frac{B_0B_1}t+{\rm tr}\left(B_1\Bigl(\begin{array}{cc} 1 & 0 \\ 0 & 0 \end{array}\Bigr)\right)+{\rm tr}\frac{B_1^2}{t+1}.
\end{eqnarray*}
Therefore,
$$
d\ln\tilde\tau(s(t))=\frac12{\rm res}_{z=t}{\rm tr}\Bigl(B(z,t)\Bigr)^2dt+{\rm tr}\frac{B_1^2}{t+1}\,dt.
$$
Denoting now $B_1(t)=(b_{ij}(t))$ we have
\begin{eqnarray*}
{\rm tr}\,B_1^2&=&(b_{11})^2+2b_{12}b_{21}+(b_{22})^2=(b_{11}+b_{22})^2-2(b_{11}b_{22}-b_{12}b_{21})=\\
&=&({\rm tr}\,B_1)^2-2\det B_1=\theta_1^2/2={\rm const}
\end{eqnarray*}
(recall that the eigenvalues of the matrix $B_1(t)$ are $\pm\theta_1/2$ not depending on $t$). Hence,
$$
d\ln\frac{\tilde\tau(s(t))}{(t+1)^{\theta_1^2/2}}=\frac12{\rm res}_{z=t}{\rm tr}\Bigl(B(z,t)\Bigr)^2dt.
$$
Thus, near a point $t^*\in\Theta$ we have the formula
\begin{eqnarray}\label{localmiwa}
d\ln\tau^*(t)=\frac12{\rm res}_{z=t}{\rm tr}\Bigl(B(z,t)\Bigr)^2dt
\end{eqnarray}
for the local $\tau$-function $\tau^*(t)=\tilde\tau(s(t))/(t+1)^{\theta_1^2/2}$ of (\ref{isomfam}).

Further we use standard reasonings to come from the local $\tau$-functions to a global one and finish the proof of the
theorem. Consider a covering $\{U_{\alpha}\}$ of the deformation space
$\mathcal{D}$ such that for every $U_{\alpha}$ there is a
holomorphic function $\tau_{\alpha}(t)$ satisfying the equality (\ref{localmiwa})
and non-vanishing in $U_{\alpha}$ if $U_{\alpha}\cap\Theta=\varnothing$. In every
non-empty intersection $U_{\alpha}\cap U_{\beta}$ one has 
\begin{eqnarray}\label{cocycle1}
\tau_{\alpha}(t)=c_{\alpha\beta}\tau_{\beta}(t),
\end{eqnarray}
where $c_{\alpha\beta}={\rm const}\ne0$, since $d\ln c_{\alpha\beta}=d\ln\tau_{\alpha}-d\ln\tau_{\beta}=0$. 
The equalities (\ref{cocycle1}) imply
\begin{eqnarray}\label{cocycle2}
c_{\alpha\beta}c_{\beta\gamma}=c_{\alpha\gamma}
\end{eqnarray}
for non-empty intersections $U_{\alpha}\cap U_{\beta}\cap U_{\gamma}$. 

Let us fix logarithms $l_{\alpha\beta}=\ln c_{\alpha\beta}$ in such a way that $l_{\alpha\beta}=-l_{\beta\alpha}$.
Then, as follows from (\ref{cocycle2}),
$$
l_{\alpha\beta}-l_{\alpha\gamma}+l_{\beta\gamma}=2\pi{\bf i}\,l_{\alpha\beta\gamma},
$$
where the set of numbers $l_{\alpha\beta\gamma}\in\mathbb Z$ defines an element of the \v{C}ech cohomology group $H^2({\mathcal D}, {\mathbb Z})$.
Since $\cal D\cong\mathbb C$, one has $H^2({\mathcal D}, {\mathbb Z})=0$, hence there is a set $\{l'_{\alpha\beta}\}\subset\mathbb Z$
such that 
$$
l_{\alpha\beta\gamma}=l'_{\alpha\beta}-l'_{\alpha\gamma}+l'_{\beta\gamma}\quad\mbox{for}\quad 
U_{\alpha}\cap U_{\beta}\cap U_{\gamma}\ne\varnothing.
$$
Therefore
$$
(l_{\alpha\beta}-2\pi{\bf i}\,l'_{\alpha\beta})-(l_{\alpha\gamma}-2\pi{\bf i}\,l'_{\alpha\gamma})
+(l_{\beta\gamma}-2\pi{\bf i}\,l'_{\beta\gamma})=0,
$$
and the set of numbers $\lambda_{\alpha\beta}=l_{\alpha\beta}-2\pi{\bf i}\,l'_{\alpha\beta}\in\mathbb C$ defines an element of the 
\v{C}ech cohomology group $H^1({\mathcal D}, {\mathbb C})$. Since the latter is trivial, there is a set 
$\{\lambda_{\alpha}\}\subset\mathbb C$ such that 
$$
\lambda_{\alpha\beta}=\lambda_{\alpha}-\lambda_{\beta}\quad\mbox{for}\quad U_{\alpha}\cap U_{\beta}\ne\varnothing.
$$
Thus, a set $\{c_{\alpha}=e^{\lambda_{\alpha}}\}$ of non-zero constants is such that
$$
c_{\alpha\beta}=e^{l_{\alpha\beta}}=e^{\lambda_{\alpha\beta}}=c_{\alpha}c_{\beta}^{-1}
$$ 
in every non-empty intersection $U_{\alpha}\cap U_{\beta}$. Therefore we have a global function $\tau(t)$ holomorphic on $\cal D$ and equal to 
$c_{\alpha}^{-1}\tau_{\alpha}(t)$ in every $U_{\alpha}$, hence satisfying (\ref{localmiwa}).

\section{Proof of Lemma 1}

Consider a point $s^*\in s(\Theta)$. Though the family (\ref{isomfam2}) is not defined for this value of the parameter, 
one can construct an auxiliary linear meromorphic $(2\times2)$-system
\begin{eqnarray}\label{linsyst_aux}
\frac{dw}{d\xi}=A^*(\xi)\,w,\qquad
A^*(\xi)=\frac{A_{01}^*}{\xi}+\frac{A_{02}^*}{\xi^2}+\frac{A_1^*}{\xi-1}+\frac{A_2^*}{\xi-s^*},
\end{eqnarray}
with irregular non-resonant singular point $\xi=0$ of Poincar\'e rank $1$ and Fuchsian singular points
$1, s^*$. In a neighbourhood of zero this system is holomorphically equivalent to the image of the initial system 
(\ref{22syst}) under the change of variable $\xi=1/(z+1)$, has the same monodromy matrices as the latter,
but it has an {\it apparent} Fuchsian singularity at the infinity ({\it i.~e.}, the monodromy at this point is trivial). 

We will use the following facts explained in \cite{BG}:
\begin{itemize}

\item the auxiliary system $(\ref{linsyst_aux})$ is included into a
(Malgrange) isomonodromic family
\begin{eqnarray}\label{isom_fam_aux}
\frac{dw}{d\xi}=\biggl(
\frac{A_{01}(s)}{\xi}+\frac{A_{02}(s)}{\xi^2}+\frac{A_1(s)}{\xi-1}+\frac{A_2(s)}{\xi-s}\biggr) w
\end{eqnarray}
whose isomonodromic fundamental matrix $W(\xi,s)$ near the point $\xi=\infty$ has the form
\begin{eqnarray}\label{fundinf}
W(\xi,s)=U(\xi,s)\xi^K, \qquad U(\xi,s)=I+U_1(s)\frac1{\xi}+U_2(s)\frac1{\xi^2}+\ldots, 
\end{eqnarray}
where $K={\rm diag}\,(-1,1)$ and $\frac{dU_1(s)}{ds}=-A_2(s)$; 

\item the upper right element $u_1(s)$ of the matrix $U_1(s)$ vanishes at the point $s=s^*$ and is
not equal to zero identically, hence it is a local $\tau$-function of the family (\ref{isomfam2}) 
(and $\frac{du_1}{ds}(s^*)\ne0$ if the monodromy of (\ref{isomfam2}) is irreducible);

\item for $s\ne s^*$ the family (\ref{isomfam2}) is meromorphically equivalent to (\ref{isom_fam_aux}) {\it via} 
a gauge transformation $\tilde y=\Gamma_1(\xi,s)w$, where
$$
\Gamma_1(\xi,s)=\left(\begin{array}{c c} \frac{f(s)}{u_1(s)}+\xi & -u_1(s) \\
                                         \frac1{u_1(s)} & 0 
											\end{array}\right),
$$
for some function $f$ holomorphic at the point $s=s^*$.
\end{itemize}

The coefficient matrices $\widetilde B(\xi,s)$ and $A(\xi,s)$ of the families (\ref{isomfam2}) and (\ref{isom_fam_aux})
respectively are connected by the equality
$$
\widetilde B(\xi,s)= \frac{\partial \Gamma_1}{\partial\xi}\Gamma_1^{-1}+
\Gamma_1 A(\xi,s) \Gamma_1^{-1},
$$
therefore
\begin{eqnarray*}
\Bigl(\widetilde B(\xi,s)\Bigr)^2&=&\Bigl(\frac{\partial\Gamma_1}{\partial\xi}
\Gamma_1^{-1}\Bigr)^2+\frac{\partial\Gamma_1}{\partial\xi}\,A(\xi,s)\Gamma_1^{-1}+\\
 & &+\Gamma_1\,A(\xi,s)\Gamma_1^{-1}\frac{\partial\Gamma_1}{\partial\xi}\Gamma_1^{-1}+
\Gamma_1\Bigl(A(\xi,s)\Bigr)^2\Gamma_1^{-1}.
\end{eqnarray*}
As follows from the form of the matrix $\Gamma_1(\xi,s)$, the
product $(\frac{\partial}{\partial\xi}\Gamma_1)\Gamma_1^{-1}$
is of the form
$$
\left(\begin{array}{c c} 1 & 0 \\
                         0 & 0 
			\end{array}\right)
\left(\begin{array}{c c} 0 & u_1 \\
                         -\frac1{u_1} & * 
			\end{array}\right)=\left(\begin{array}{c c} 0 & u_1 \\
                                                  0 & 0 
			                         \end{array}\right),			
$$
hence its square is the zero matrix and
$$
{\rm tr}\Bigl(\widetilde B(\xi,s)\Bigr)^2=2\,{\rm tr}\Bigl(
\Gamma_1^{-1}\frac{\partial\Gamma_1}{\partial\xi}\,A(\xi,s)\Bigr)
+{\rm tr}\Bigl(A(\xi,s)\Bigr)^2.
$$
Since 
$$
\Gamma_1^{-1}\frac{\partial\Gamma_1}{\partial\xi}=
\left(\begin{array}{c c} 0 & 0 \\
                         -\frac1{u_1} & 0 
			\end{array}\right),
$$
one has
$$
{\rm tr}\Bigl(\widetilde B(\xi,s)\Bigr)^2=-2\,\frac{a(\xi,s)}{u_1(s)}
+{\rm tr}\Bigl(A(\xi,s)\Bigr)^2,
$$
where $a(\xi,s)$ is the upper right element of the matrix $A(\xi,s)$. Thus,
$$
{\rm res}_{\xi=s}{\rm tr}\Bigl(\widetilde B(\xi,s)\Bigr)^2ds=-2\,\frac{{\rm res}_{\xi=s}a(\xi,s)}{u_1(s)}\,ds+
{\rm res}_{\xi=s}{\rm tr}\Bigl(A(\xi,s)\Bigr)^2ds.
$$
As ${\rm res}_{\xi=s}a(\xi,s)$ is equal to the upper right element of $A_2(s)$, which is $-\frac{du_1}{ds}$, one has 
$$
{\rm res}_{\xi=s}{\rm tr}\Bigl(\widetilde B(\xi,s)\Bigr)^2ds=2\,d\ln u_1(s)+
{\rm res}_{\xi=s}{\rm tr}\Bigl(A(\xi,s)\Bigr)^2ds.
$$

The differential $1$-form ${\rm res}_{\xi=s}{\rm tr}\Bigl(A(\xi,s)\Bigr)^2ds$
is closed and holomorphic in a neighbourhood $D(s^*)$ of the point $s=s^*$, hence
there is a function $f^*$ holomorphic and non-vanishing in $D(s^*)$ such that 
$d\ln f^*(s)=\frac12{\rm res}_{\xi=s}{\rm tr}\Bigl(A(\xi,s)\Bigr)^2ds$. Therefore,
$$
{\rm res}_{\xi=s}{\rm tr}\Bigl(\widetilde B(\xi,s)\Bigr)^2ds
=2\,d\ln (f^*(s)u_1(s)),
$$
which finishes the proof of the lemma, with $\tilde\tau(s)=f^*(s)u_1(s)$.
\bigskip

{\bf Acknowledgements.} This work is supported by the Russian Foundation for Basic Research (grant no. RFBR-14-01-31145 mol\_a) 
and RF President program for young scientists (grant no. MK-4594.2013.1).

\end{document}